\documentclass[11pt,reqno]{amsart}
\usepackage[utf8]{inputenc}
\usepackage{amsmath,amsthm,amssymb,tikz-cd,mathrsfs,enumitem,caption,subcaption}
\usepackage[scale=0.75]{geometry}

\title[Misiurewicz poly. for rational maps with nontrivial aut. II]{Misiurewicz polynomials for rational maps with nontrivial automorphisms II}
\author{Minsik Han}
\address{Department of Mathematics, Box 1917, Brown University, Providence, RI 02912, USA}
\email{minsik\_han@brown.edu}

\makeatletter
\@namedef{subjclassname@2020}{%
  \textup{2020} Mathematics Subject Classification}
\makeatother
\subjclass[2020]{37P05}
\keywords{Arithmetic dynamics, Rational maps, Automorphisms, Gleason polynomials, Misiurewicz polynomials, Irreducibility}

\newcommand{\bP}{\mathbb{P}}
\newcommand{\bQ}{\mathbb{Q}}
\newcommand{\bZ}{\mathbb{Z}}

\newcommand{\cF}{\mathcal{F}}
\newcommand{\cG}{\mathcal{G}}
\newcommand{\cL}{\mathcal{L}}
\newcommand{\cM}{\mathcal{M}}
\newcommand{\cS}{\mathcal{S}}

\newcommand{\ord}{\textup{ord}}

\newtheorem{theorem}{Theorem}

\newtheorem{lemma}[theorem]{Lemma}

\newtheorem{proposition}[theorem]{Proposition}
\newtheorem{maintheorem}{Theorem}

\theoremstyle{remark}
\newtheorem*{remark}{Remark}

\begin{document}

\begin{abstract}
    This paper continues discussions in the author's previous paper about the Misiurewicz polynomials defined for a family of degree $d \ge 2$ rational maps with an automorphism group containing the cyclic group of order $d$. In particular, we extend the sufficient conditions that the Misiurewicz polynomials are irreducible over $\bQ$. We also prove that the Misiurewicz polynomials always have an irreducible factor of large degree.
\end{abstract}

\maketitle

\section{Introduction and main results}

A rational map on $\bP^1$ is called post-critically finite if all of its critical points are preperiodic, i.e., they have finite forward orbits. For a family of unicritical polynomial maps $z \mapsto z^d+c$, the $c$-values which make the map post-critically finite are called the Misiurewicz points. The Misiurewicz points for a fixed degree $d$ and a fixed dynamical portrait $(m,n)$ are roots of a polynomial called the Misiurewicz polynomial. The irreducibility of Misiurewicz polynomials for unicritical polynomial maps are widely studied \cite{buff,goksel,goksel2} but still remains open in general.

The author discussed a similar problem in \cite{han}, but instead for a family of rational maps. Taking rational maps with a specific automorphism group, it was possible to control the dynamic aspects of their critical points. Explicitly, in the moduli space $\cM_d$ of degree $d$ rational maps on $\bP^1$, a single-parameter family \[\left\{\phi_a = \frac{az}{z^d+d-1} : a \ne 0 \right\}\] parametrizes the maps with an automorphism group containing the cyclic group of order $d$. We constructed a polynomial $G_m := G_{m,1}\in \bZ[a]$ whose roots make all nontrivial critical points of $\phi_a$ post-critically finite with the dynamical portrait $(m,1)$.

In \cite{han}, we proved that the polynomial $G_m$ is irreducible if $d$ is an odd prime and $m$ is at most $3$. We extend this result in this paper to all $m\le d$, and prove a general statement for larger $m$. Here is our main theorem.

\begin{maintheorem}
\label{mainthm}
Fix a prime $d \ge 3$. Then there exists an irreducible factor $F_m$ of the Misiurewicz polynomial $G_m:=G_{m,1}$ corresponding to the family $\{\phi_a\}$ over $\bQ$, where \[ \liminf_{m\to \infty} \frac{\deg F_m}{\deg G_m} \ge 1-\frac{1}{d^d}.\] Moreover, if $m\le d$, then $G_m$ is irreducible over $\bQ$.
\end{maintheorem}

\begin{remark}
The above statement is also true over $\bQ_d$. It has been conjectured that the Misiurewicz polynomial is always irreducible over $\bQ$, while the polynomial is not irreducible over $\bQ_d$ in general. For example, if $d=3$ and $m=4$ then $G_4$ is a $55$-degree polynomial which is irreducible over $\bQ$, but it is factored as a product of one $53$-degree polynomial and two linear polynomials over $\bQ_3$.
\end{remark}

\section{Backgrounds}

The map $\phi_a$ is represented with homogeneous coordinates as \[\phi_a ([x,y]) = [axy^{d-1},x^d+(d-1)y^d].\]
We first do a change of variable $a\mapsto (b+1)d$, so we work on
\[\varphi_b([x,y]) = [(b+1)dxy^{d-1},x^d+(d-1)y^d].\]
Then the Misiurewicz polynomials for the dynamical portrait $(m,1)$ is defined as \[\cG_m= \frac{1}{bd} \left[(b+1)^d d^d s_{m-1}^{d(d-1)} - (bd+1) \frac{s_m^d-((b+1)d s_{m-1}^d)^d}{s_m-(b+1)d s_{m-1}^d}\right] \in \bZ[b],\] where \[\phi_b^n ([1,1]) = [r_n,s_n].\] Note that $r_n$ and $s_n$ are defined by the recurrence equations \begin{equation} \label{eq:rsrecur} \begin{split} r_{n+1} &= (b+1) d r_n s_n^{d-1}, \\ s_{n+1} & = r_n^d + (d-1) s_n^d.\end{split} \end{equation}

\begin{remark}
See \cite{han} for the detailed definition of the Misiurewicz polynomials, as well as the reason why the change of variable is required.
\end{remark}

We first determine the degrees of $r_n$ and $s_n$, which will be used to prove our main result.

\begin{proposition}
\label{rsdeg}
Let $n \ge 2$ and $i,j$ be unique integers satisfying
\begin{equation} \label{nij} n=di+j,\ \ j \in \{2,\cdots,d+1\}.\end{equation}
Then \begin{equation} \label{eq:rsdeg} \deg r_n = \frac{d^{n+d-1}-d^{j-1}}{d^d-1} - (d-j),\ \ \deg s_n = \frac{d^{n+d-1}-d^{j-1}}{d^d-1}.\end{equation}
\begin{proof}
Note that by definition and the fact that the leading coefficients of $r_n$ and $s_n$ are always positive, we have \[\deg r_n = 1+\deg r_{n-1} + (d-1) \deg s_{n-1} ,\ \ \deg s_n = d\max (\deg r_{n-1},\deg s_{n-1}).\] It follows that
\begin{equation} \label{eq:rsdegrecur} \begin{split} \deg s_{n-1} \ge \deg r_{n-1} \ \ \Longrightarrow \ \ &\deg s_n = d \deg s_{n-1}, \\ &\deg s_n - \deg r_n = \deg s_{n-1} - \deg r_{n-1}-1. \end{split} \end{equation}

We prove by induction on $i$. If $i=0$, we need to show \begin{equation} \label{eq:rsi0} \deg r_j = d^{j-1} - (d-j),\ \ \deg s_j = d^{j-1}\end{equation} for all $j$. We can prove the $j=2$ case directly, \[\deg r_2 = 2,\ \ \deg s_2 = d.\] Then applying \eqref{eq:rsdegrecur} repeatedly, we get \eqref{eq:rsi0} for all $j$, up to \[\deg r_{d+1} = d^d+1,\ \ \deg s_{d+1} = d^d.\]

Now assume that we proved \eqref{eq:rsdeg} up to $i$. In particular, \[\deg r_{di+(d+1)} = \frac{d^{di+2d}-d^d}{d^d-1}+1,\ \ \deg s_{di+(d+1)} =  \frac{d^{di+2d}-d^d}{d^d-1}.\] Then the next terms are given by \begin{align*} \deg r_{d(i+1)+2} & = 2 + d \deg s_{di+(d+1)} \\ &= 2 +  \frac{d^{di+2d+1}-d^{d+1}}{d^d-1} \\&= \frac{d^{(d(i+1)+2)+d-1} - d}{d^d-1} -(d-2)\end{align*} and \begin{align*} \deg s_{d(i+1)+2} & = d \deg r_{di+(d+1)} \\ &= \frac{d^{di+2d+1}-d^{d+1}}{d^d-1} + d \\ &= \frac{d^{(d(i+1)+2)+d-1} - d}{d^d-1},\end{align*} which coincide with \eqref{eq:rsdeg}. Moreover, using \eqref{eq:rsdegrecur} again we can conclude that \eqref{eq:rsdeg} holds for $i+1$ and all $j$ as well. This ends the proof.
\end{proof}
\end{proposition}

\section{Newton Polygons and Notations}

For a prime $p$ and a polynomial $f(z)=\sum_{i=0}^n a_i z^i \in \bQ [z]$ where $a_n \ne 0$, the \emph{$p$-Newton polygon} of $f$, denoted by $N_p(f)$, is defined as the lower convex hull of the set of points \[\{(i,\ord_p(a_i)):i=0,1,\cdots,n\},\] where $\ord_p$ is the $p$-adic valuation on $\bQ$. By definition, $N_p(f)$ is a union of connected line segments with increasing slopes. Note that the slope may be $-\infty$, in case that $a_0 = 0$. We use a notation \[v_{i,p}(f):=\ord_p(a_i),\] which is often simplified as $v_i(f)$ if $p$ is obvious. Also, we define \[V_p(f):= \min \{v_{i,p}(f): i=0,1,\cdots,n\},\] which is simplified as $V(f)$ as well. Note that \[V(fg) = V(f)V(g),\ \ V(f+g) \ge \min(V(f),V(g))\] for any pair of polynomials $(f,g)$.

We introduce some more notation which simplifies our statements. The points on $N_p(f)$ where the slope changes are called \emph{vertices} of $N_p(f)$. Two endpoints $(0,v_0(f))$ and $(n,v_n(f))$ are also considered as vertices. The \emph{principal $p$-Newton polygon}, denoted by $N_p^-(f)$, is the subpolygon of $N_p(f)$ which consists of line segments with negative slopes, including $-\infty$. In this paper, we represent $N_p^-(f)$ by vertices of it with finite coordinates as \[N_p^-(f) = L\Bigl(\bigl(i_1,v_{i_1}(f)\bigr),\bigl(i_2,v_{i_2}(f)\bigr),\cdots,\bigl(i_k,v_{i_k}(f)\bigr)\Bigr),\] where \[i_1 = \min\{i:v_i(f)<\infty\}\] and $i_1<i_2<\cdots<i_k$. Here $(i_1,v_{i_1}(f))$ is called the \emph{initial point} of $N_p^-(f)$.

\begin{figure}[ht]
    \centering
    \begin{subfigure}[h]{0.495\textwidth}
        \centering
        \includegraphics[height=4.5cm]{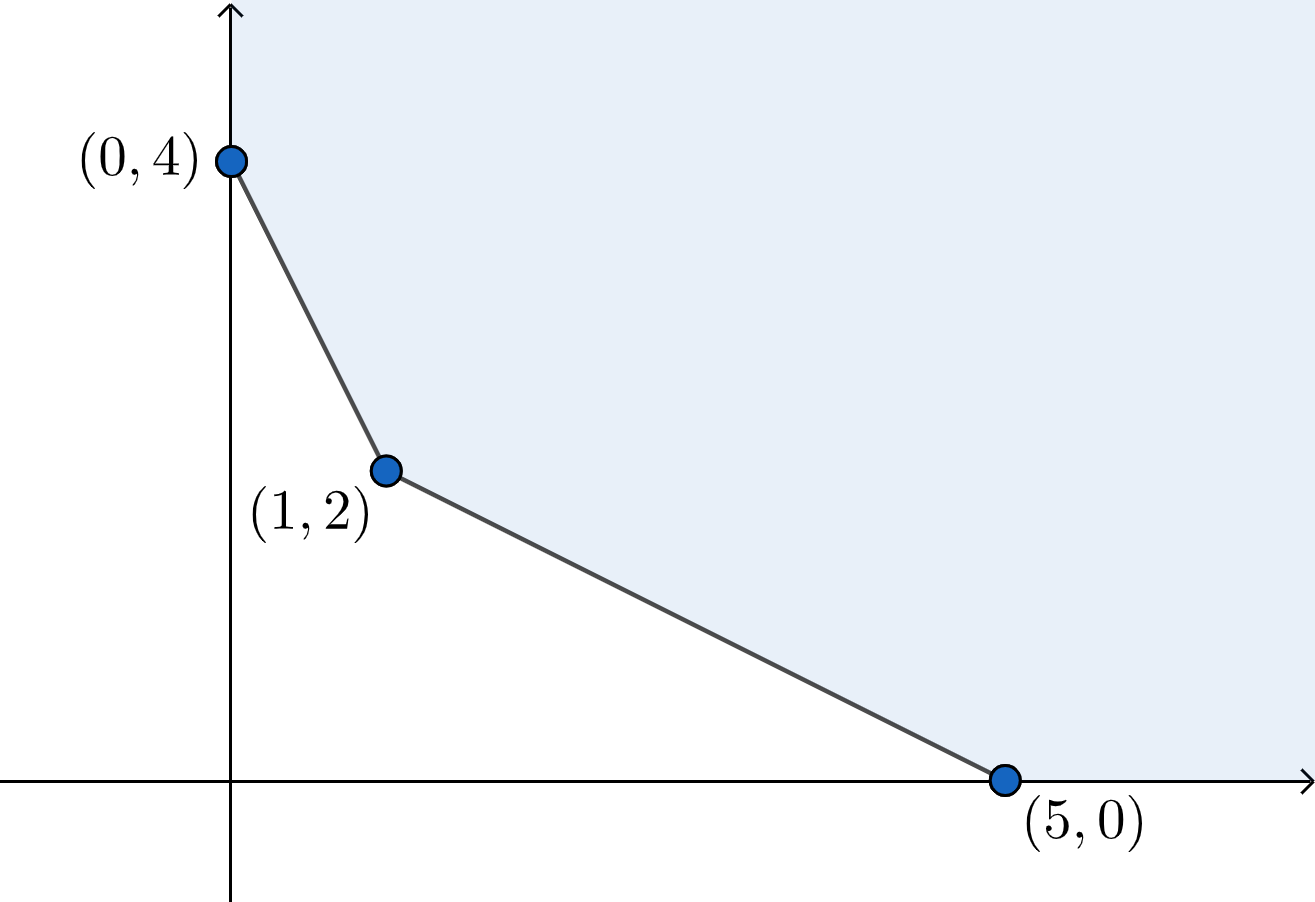}
        \caption{$N_p^-(f_1) = L\bigl((0,4),(1,2),(5,0)\bigr)$}
    \end{subfigure}
    \hfill
    \begin{subfigure}[h]{0.495\textwidth}
        \centering
        \includegraphics[height=4.5cm]{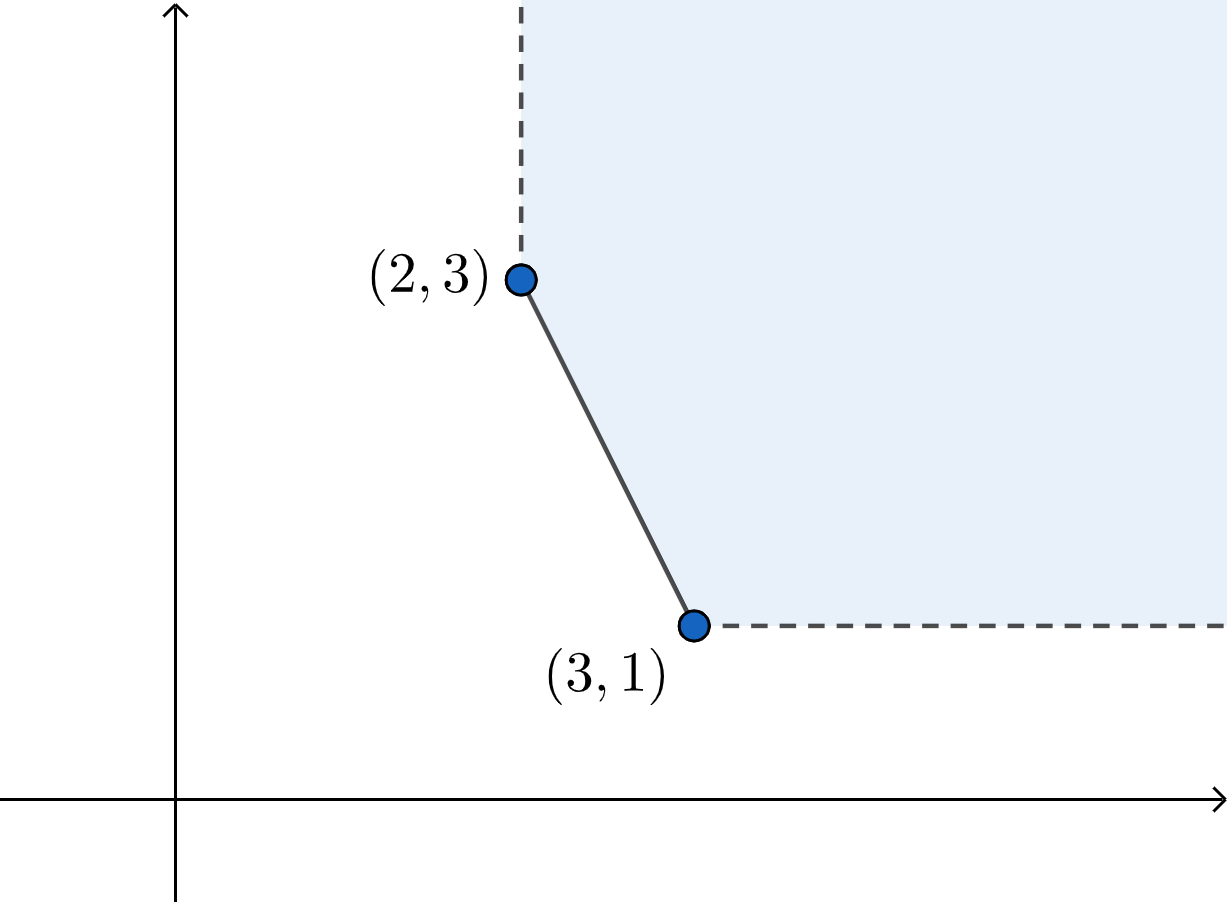}
        \caption{$N_p^-(f_2) = L\bigl((2,3),(3,1)\bigr)$}
    \end{subfigure}
    \caption{Principal $p$-Newton polygons}
    \label{fig:convexhull}
\end{figure}

For example, for $p=2$ and two polynomials \[f_1(z) = z^5 + 4z^3 + 4z-16,\ \ f_2(z) = 6z^4-2z^3+8z^2\] we have \[N_p^-(f_1) = L\bigl((0,4),(1,2),(5,0)\bigr),\ \ N_p^-(f_2) = L\bigl((2,3),(3,1)\bigr).\] See Figure \ref{fig:convexhull} below for illustrations.

For given finite number of principal $p$-Newton polygons, we can define the sum of them as follows: the sum of all initial points is the initial point of the sum, and from there we attach all line segments from all polygons in the order of increasing slopes. For any polynomials $f_1,\cdots,f_l$, the principal $p$-Newton polygon of the product is equal to the sum of principal $p$-Newton polygons, that is, \begin{equation} \label{eq:polygonsum} N_p^-(f_1 \cdots f_l) = N_p^-(f_1) + \cdots + N_p^-(f_l).\end{equation} For a detailed study about principal Newton polygons, see \cite{guardia}. Note that \eqref{eq:polygonsum} also holds for the original Newton polygons as well. \cite[Proposition 5]{han}

We finish this section with a useful lemma, which was proved and used in \cite{han}.

\begin{lemma}
\label{orderlemma}
Let $f(z)=a_n z^n + \cdots + a_0$ be a polynomial, and consider the $k$\textsuperscript{th} power $f(z)^k$. Let $(\alpha_1,\cdots,\alpha_k)$ be a $k$-tuple of integers satisfying $0\le \alpha_1 \le \cdots \le \alpha_k \le n$, and let $N(\alpha_1,\cdots,\alpha_k)$ be the number of permutations of the subscripts maintaining the ordering. That is, $N(\alpha_1,\cdots,\alpha_k)$ is the number of permutations $(r_1,\cdots,r_k)$ of $(1,\cdots,k)$ such that $\alpha_{r_1} \le \cdots\le \alpha_{r_k}$. Then
\begin{equation} \label{eq13} v_i(f(z)^k)\ge \min_{\substack{0 \le \alpha_1 \le \cdots \le \alpha_k \\ \alpha_1+\cdots+\alpha_k = i}}\left( \ord_p \left (\frac{k!}{N(\alpha_1,\cdots,\alpha_k)}\right)+\sum_{j=1}^k v_{\alpha_j}(f)\right).\end{equation}
\end{lemma}

In particular, if $k=p$ and $i$ is not divisible by $p$ then \begin{equation} \label{eq:ordernmid} v_i(f(z)^p) \ge 1+ \min_{\substack{0 \le \alpha_1 \le \cdots \le \alpha_p \\ \alpha_1+\cdots+\alpha_p = i}} \left( \sum_{j=1}^p v_{\alpha_j}(f)\right),\end{equation} and if $i=pe$ is divisible by $p$ then \begin{equation} \label{eq:ordermid} v_i(f(z)^p) \ge \min \left[ 1+ \min_{\substack{0 \le \alpha_1 \le \cdots \le \alpha_p \\ \alpha_1+\cdots+\alpha_p = i}} \left( \sum_{j=1}^p v_{\alpha_j}(f)\right), p v_e (f)\right].\end{equation}

\ifx
For any $\lambda=-a/b \in \bQ^-$ such that $a,b$ are coprime, let $\cS(\lambda)$ be the set of all line segments with slope $\lambda$ whose endpoints are lattice points with nonnegative coordinates. We also include lattice points with nonnegative coordinates in $\cS(\lambda)$. For any $S \in \cS(\lambda)$, let $i(S)$ be the initial point, i.e., the upper left endpoint of $S$, and $l(S),h(S)$ be the horizontal length and height of $S$, respectively. Then we can define the degree of $S$ \[d(S) = \frac{l(S)}{b} = \frac{h(S)}{a}.\] Note that for any pair $(i,d)$ there exists at most one $S$ such that $i(S)=i$ and $d(S)=d$. Then we define the addition on $\cS(\lambda)$ in a natural way, \[(i,d)+(i',d') = (i+i',d+d').\]

Next, we define the set $\cS(-\infty)$ of all line segment with slope $-\infty$. An element $S$ of $\cS(-\infty)$, roughly an upward ray from a point, is determined by the initial point $i(S)$ and a formal horizontal length $l(S)$. Then we can similarly define the addition on $\cS(-\infty)$.

Now let $\cS = \cS(-\infty) \sqcup \left( \bigcup_{\lambda \in \bQ^-} \cS(\lambda)\right)$ and consider the formal sums on $\cS$.

Let \[S=\{(x_i,y_i): 0\le x_1 < \cdots < x_m\} \subset \bZ^2 \] be a set of lattice points on the coordinate plane, indexed by $x$-coordinates. We say that $S$ is \emph{de-convex} if the slopes of line segments between two adjacent points are increasing and negative, i.e., \[\frac{y_2-y_1}{x_2-x_1} \le \cdots \le \frac{y_m-y_{m-1}}{x_m-x_{m-1}}<0.\] For a de-convex set of lattice points $S$, $L(S)=L\bigl((x_1,y_1),\cdots,(x_m,y_m)\bigr)$ denotes the lower convex hull of $S \cup \{(\infty,y_m)\}$.

Let $f(z)=\sum_{n=0}^N a_n z^n \in \bQ[z]$ be a polynomial and $S=\{(x_1,y_1),\cdots,(x_m,y_m)\}$ be a de-convex set of lattice points. Then we say that $N_p(f) \ge L(S)$ if $a_n=0$ for all $n<x_1$, and the $p$-Newton polygon of $f$ is contained in the closed convex region defined by $L(S)$. In other words, $N_p(f) \ge L(S)$ if and only if $v_{i,p}(f) \ge \cL(i)$ for all $i\ge 0$, where $\cL$ is a piecewise function on $[0,\infty)$ defined as following.
\[\cL(x) = \begin{cases}\infty & \text{if } 0\le x<x_1, \\ \displaystyle \frac{y_{i+1}-y_i}{x_{i+1}-x_i}(x-x_i)+y_i & \text{if } x_i \le x < x_{i+1}, \\ y_m & \text{if } x_m \le x.\end{cases}\]

Now let $S,S'$ be two de-convex sets of lattice points. As in the case of Newton polygons, we can define a product $L(S)L(S')$ using the `rearranged concatenation.' Then if $f,g\in \bQ[z]$ satisfy $N_p(f) \ge L(S)$ and $N_p(g) \ge L(S')$, then $N_p(fg) \ge L(S)L(S')$.
\fi

\section{Proof of the main result}

In this section, we prove Theorem \ref{mainthm}. From the equation \[bd \cG_m = (b+1)^d d^d s_{m-1}^{d(d-1)} - (bd+1) \frac{s_m^d - ((b+1)d s_{m-1}^d)^d}{s_m-(b+1)d s_{m-1}^d},\] let \begin{equation} \label{eq:sigmatau} \sigma_m = (b+1)d s_{m-1}^d,\ \ \tau_m= s_m-\sigma_m.\end{equation} Then from above,
\begin{equation} \begin{aligned} \label{gsigmatau} bd \cG_m &= (b+1)d \sigma_m^{d-1} - (bd+1) \frac{s_m^d - \sigma_m^d }{s_m-\sigma_m} \\ &= (b+1)d \sigma_m^{d-1} - (bd+1) \frac{(\sigma_m+\tau_m)^d - \sigma_m^d}{\tau_m} \\&= (b+1)d \sigma_m^{d-1} - (bd+1) \sum_{k=0}^{d-1} \binom{d}{k} \sigma_m^k \tau_m^{d-1-k} \\ &= -bd(d-1)\sigma_m^{d-1} - (bd+1) \sum_{k=0}^{d-2} \binom{d}{k} \sigma_m^k \tau_m^{d-1-k} \\ &= -\sum_{k=0}^{d-1} \cF_k, \end{aligned} \end{equation}
where \[\cF_k = \begin{cases} \displaystyle (bd+1) \binom{d}{k} \sigma_m^k \tau_m^{d-1-k} & \text{if } k=0,\cdots,d-2,\\ bd(d-1) \sigma_m^{d-1} & \text{if } k=d-1.\end{cases}\]

\begin{proposition}
\label{sigmataudegree}
If $m \equiv 2 \pmod d$, then \[\deg \sigma_m = \frac{d^{m+d-1}-d^{d+1}}{d^d-1}+1,\ \ \deg \tau_m = \frac{d^{m+d-1}-d}{d^d-1} = \deg s_m.\] In other cases, \[\deg \sigma_m = \deg \tau_m = \frac{d^{m+d-1}-d^{j-1}}{d^d-1}+1,\] where $j$ is defined in \eqref{nij}. In particular, \[\lim_{m\to \infty} \frac{\deg \sigma_m}{d^{m-1}}= \lim_{m\to \infty} \frac{\deg \tau_m}{d^{m-1}} = \frac{d^d}{d^d-1}.\]
\begin{proof}
If $m \not \equiv 2 \pmod d$, then $j>2$ so from \eqref{eq:sigmatau}, \[\deg \sigma_m = d \deg s_{m-1} +1 = \frac{d(d^{m+d-2}-d^{j-2})}{d^d-1} +1 = \frac{d^{m+d-1}-d^{j-1}}{d^d-1}+1 = \deg s_m +1,\] thus $\deg \tau_m = \deg \sigma_m$.

On the other hand, if $m \equiv 2 \pmod d$, then $m=di+2$ for some $i$, so \[m-1 = d(i-1)+(d+1).\] Therefore, \[\deg \sigma_m = d \deg s_{m-1} +1 = \frac{d(d^{m+d-2} - d^d)}{d^d-1} +1 = \frac{d^{m+d-1} - d^{d+1}}{d^d-1} +1  = \deg s_m - (d-1),\] so $\deg \tau_m = \deg s_m$.
\end{proof}
\end{proposition}

Now we investigate the principal $d$-Newton polygon of these polynomials. Let \[D_{m} = d^{m-1}+ d^{m-2}+\cdots+1 = \frac{d^m-1}{d-1}.\] Note that $dD_{m-1}+1 = D_m$, so \[\frac{D_{m-1}}{d^{m-1}} = \frac{dD_{m-1}}{d^m} < \frac{D_m}{d^m}\] for every $m$.

\begin{proposition}
\label{smnp}
Let $m\ge 2$. Then \begin{equation} \label{eq:smnp} N_d^-(s_m) = L\bigl((0,D_m),(d^{m-1},d^{m-1})\bigr).\end{equation} That is, $N_d^-(s_m)$ has the initial point $(0,D_m)$ and it consists of a single line segment with slope $-D_{m-1}/d^{m-1}$.
\begin{proof}
We prove by induction. First, \[s_2 = ((b+1)d)^d + (d-1)d^d = ((b+1)^d + (d-1))d^d.\] From \[(b+1)^d + (d-1) = b^d + \binom{d}{d-1} b^{d-1} + \cdots + \binom{d}{1} b + d,\] it follows that $N_d(s_2)$ has exactly two vertices $(0,d+1)$ and $(d,d)$ so \[N_d^-(s_2)= L\bigl((0,D_2),(d,d)\bigr).\]

Now assume that \eqref{eq:smnp} holds for all values up to $m$. We show that the same holds for $m+1$, investigating $v_i(s_{m+1})$ for values of $i$. Explicitly, we need to show that \[v_i(s_{m+1}) \begin{cases} = D_{m+1} & \text{if } i=0, \\ \ge \displaystyle D_{m+1} - \frac{D_m i}{d^m} & \text{if } 0<i<d^m, \\ = d^m & \text{if }i=d^m, \\ \ge d^m & \text{if }d^m<i.\end{cases}\]

\begin{enumerate}[label=\textbf{(\Alph*)}]
    \item $i=0$: We observe that if $r_m \equiv s_m \pmod b$, then \[r_{m+1} \equiv s_{m+1} \equiv d s_m^d \pmod b\] from \eqref{eq:rsrecur}. Since $r_0 = s_0 = 1$, it follows that \[v_0(s_{m+1}) = 1+ d v_0 (s_m) = 1+dD_m = D_{m+1}.\] This proves the first case.
    
    \item $0<i<d^m$: We first claim that \[v_i (s_m^d) \ge D_{m+1} - \frac{D_m i}{d^m}.\] Note that the assumption assures that \[v_j (s_m) \ge D_m - \frac{D_{m-1}j}{d^{m-1}}\] for all $j$. If $i$ is not divisible by $d$, then from \eqref{eq:ordernmid} \[v_i (s_m^d) \ge 1+ dD_m - \frac{D_{m-1}i}{d^{m-1}} \ge D_{m+1} - \frac{D_m i}{d^m}.\]
    
    Therefore, it suffices to show that \[dv_e(s_m) \ge  D_{m+1} -  \frac{D_mde}{d^m} = D_{m+1} - \frac{D_m e}{d^{m-1}} = D_{m+1} - e- \frac{D_{m-1}e}{d^{m-1}}\] for all $e=1,\cdots,d^{m-1}-1$. We actually prove a stronger inequality, \[dv_e(s_m) \ge D_{m+1}-e - \left \lfloor \frac{D_{m-1}e}{d^{m-1}}\right \rfloor.\] Again by the assumption, we have \[v_e(s_m) \ge D_m - \frac{D_{m-1}e}{d^{m-1}},\] and since $v_e(s_m)$ is an integer, \[v_e(s_m) \ge D_m - \left\lfloor \frac{D_{m-1}e}{d^{m-1}} \right \rfloor .\] Then we only need to show that \[d\left(D_m - \left\lfloor \frac{D_{m-1}e}{d^{m-1}} \right \rfloor \right) \ge D_{m+1} - e- \left\lfloor \frac{D_{m-1}e}{d^{m-1}}\right\rfloor \] or equivalently
    \begin{equation}
    \label{smeq1} e\ge 1+ (d-1) \left \lfloor \frac{D_{m-1}e}{d^{m-1}}\right\rfloor.
    \end{equation}
    Now choose unique integers $u,v$ such that $e=(d-1)u+v$ and $0\le v \le d-2$. Then \[ \frac{D_{m-1}e}{d^{m-1}} = \frac{(d^{m-1}-1)u+D_{m-1}v}{d^{m-1}} = u + \frac{D_{m-1}v-u}{d^{m-1}}\] so \eqref{smeq1} is equivalent to \[v \ge 1+(d-1) \lfloor E \rfloor\] where \[E= \frac{D_{m-1}v-u}{d^{m-1}}.\] Note that \[D_{m-1}v < D_{m-1}(d-1) = d^{m-1}-1,\] so $\lfloor E\rfloor$ is at most zero. Therefore, \eqref{smeq1} naturally holds if $v\ge 1$, but it also holds if $v=0$ since then $u>0$ so $\lfloor E \rfloor$ cannot be zero. This proves the first claim.
    
    Next, we claim that \[v_i (r_m) \ge D_m - \frac{D_{m-1} i}{d^{m-1}}\] as $s_m$, thus \[v_i(r_m^d) \ge D_{m+1}-\frac{D_m i}{d^m}\] as well. From \eqref{eq:rsrecur} we have \[r_m = ((b+1)d)^m r_0 (s_0 s_1\cdots s_{m-1})^{d-1} = ((b+1)d)^m d^{d-1} (s_2 \cdots s_{m-1})^{d-1}.\] Using induction hypotheses, we get that $N_d^-(r_m)$ has the initial point \[ \left(0, (m+d-1) + (d-1) \sum_{j=2}^{m-1} D_j\right) = \left( 0, (m+d-1) + \sum_{j=2}^{m-1} (d^j-1)\right)= (0,D_m),\] which is the same as $N_d^-(s_m)$, and it consists of line segments with slopes \[-\frac{D_j}{d^j},\ \ j=1,\cdots,m-2.\] However, all this slopes are bigger than $-D_{m-1}/d^{m-1}$, which implies that $N_d^-(r_m)$ is above $N_d^-(s_m)$. This proves our next claim.
    
    Then finally, since $s_{m+1} = r_m^d + (d-1)s_m^d$, we have \[v_i(s_{m+1}) \ge \min (v_i(r_m^d),v_i(s_m^d)) \ge D_{m+1} - \frac{D_m i}{d^m}\] as desired.
    
    \item $i=d^m$: $N_d^-(s_m)=L\bigl((0,D_m),(d^{m-1},d^{m-1})\bigr)$ implies $N_d^-(s_m^d) = L\bigl((0,dD_m),(d^m,d^m)\bigr)$, so \[v_{d^m}(s_m^d) = d^m.\] On the other hand, from \eqref{eq:rsrecur} we get \[V(r_{n+1}) = 1+ V(r_n)+(d-1)V(s_n),\ \ V(s_{n+1}) \ge d \min (V(r_n),V(s_n)).\] Since $V(r_1)=V(s_1) = 1$, we can inductively show that \[V(r_m)>V(s_m) \ge d^{m-1},\] which means $V(r_m^d) > d^m$, in particular $v_{d^m}(r_m^d) > d^m$. Therefore \[v_{d^m}(s_{m+1}) = d^m.\]
    
    \item $d^m<i$: From the proof of (C) it follows that $V(s_{m+1}) \ge d^m$, which proves the last case.
\end{enumerate}
\end{proof}
\end{proposition}

It directly follows from \eqref{eq:sigmatau} that \begin{equation} \label{sigmaNewton} N_d^-(\sigma_m)= L\bigl((0,D_m),(d^{m-1},d^{m-1}+1)\bigr).\end{equation} Note that it consists of a single line segment with slope \[- \frac{dD_{m-2}}{d^{m-1}} = -\frac{D_{m-2}}{d^{m-2}}.\]

For $\tau_m$, we can observe that it is related to the Misiurewicz polynomial in a different way.

\begin{proposition}
\label{taug}
For any $m\ge 1$, \[\frac{\tau_{m}}{\tau_{m-1}} = bd \cG_{m-1},\] so \begin{equation} \tau_m = -(bd)^m \cG_{m-1} \cdots \cG_1.\end{equation}
\begin{proof}
Expressing every term using $s_i$'s, it is equivalent to \[s_m - (b+1)ds_{m-1}^d = (b+1)^d d^d s_{m-2}^{d(d-1)} (s_{m-1}-(b+1)d s_{m-2}^d) - (bd+1) (s_{m-1}^d - ((b+1)ds_{m-2}^d)^d).\] Then using $s_m = r_{m-1}^d + (d-1)s_{m-1}^d$ it is again equivalent to \begin{align*} r_{m-1}^d & = (b+1)^d d^d s_{m-2}^{d(d-1)} (s_{m-1}-(b+1)d s_{m-2}^d) + (bd+1)((b+1)ds_{m-2}^d)^d \\ &= (b+1)^d d^d  s_{m-2}^{d(d-1)} ( s_{m-1} - (d-1) s_{m-2}^d) \\ &= (b+1)^d d^d s_{m-2}^{d(d-1)} r_{m-2}^d,\end{align*} which is true by definition.
\end{proof}
\end{proposition}

\begin{theorem}
\label{gtauNewton}
For any $m \ge 1$, the followings hold.
\begin{enumerate}[label={(\roman*)}]
    \item $\cG_m$ satisfies \begin{equation} \label{gNewton} N_d^-(\cG_m)= L\bigl((0,d^m-1),(d^m-d^{m-1}-1,d^m-d^{m-1}-1)\bigr).\end{equation}
    \item $\tau_m$ satisfies \begin{equation} \begin{aligned} \label{tauNewton} N_d^-(\tau_m) & = L\bigl((d^i+(m-1-i),d^iD_{m-i}):i=0,\cdots,m-1\bigr)\\& = L\bigl((m,D_m),(d+(m-2),dD_{m-1}),\cdots,(d^{m-1},d^{m-1})\bigr).\end{aligned} \end{equation}
\end{enumerate}
\begin{proof}
We prove two statements at once by induction. If $m=1$, (i) follows from \cite{han} and (ii) follows from the fact that $\tau_1 = -bd$.

Now we assume that both statements hold for all values up to $m-1$ and prove them for $m$. Note that the principal $d$-Newton polygon given in \eqref{tauNewton} consists of line segment with slopes \[-\frac{1}{d-2}>-\frac{d}{d^2-d-1} > \cdots > -\frac{d^{m-2}}{d^{m-1}-d^{m-2}-1},\] so (ii) follows from Proposition \ref{taug} and induction hypotheses for $N_d^-(\cG_{m-1}),\cdots,N_d^-(\cG_1)$. Therefore, it suffices to prove (i), or equivalently \begin{equation} \label{eq:bdg} N_d^-(bd\cG_m) = L\bigl((1,d^m), (d^m-d^{m-1},d^m-d^{m-1})\bigr),\end{equation} with induction hypotheses and assuming (ii) for $m$ as well. Note that the principal $d$-Newton polygon given in the right hand side consists of a single line segment with slope \[-\frac{d^{m-1}}{d^m-d^{m-1}-1}.\]

Now we investigate principal $d$-Newton polygons of $\cF_k$ in \eqref{gsigmatau} for each $k$. Let $\ell$ be the line passing through two points given in \eqref{eq:bdg}. We first claim that all vertices of $N_d^-(\cF_k)$ are on or above $\ell$ for every $k$.

\begin{enumerate}[label=\textbf{(\Alph*)}]
    \item $k=0$: We have \[\cF_0 = (bd+1) \tau_m^{d-1},\] so from (ii)
    {\small \begin{align*}  N_d^-(\cF_0) & = L\Bigl(\bigl((d^i+(m-1-i))(d-1),d^iD_{m-i}(d-1)\bigr):i=0,\cdots,m-1\Bigr) \\ &= L\bigl((m(d-1),d^m-1),((d+(m-2))(d-1),d(d^{m-1}-1)),\cdots,(d^{m-1}(d-1),d^{m-1}(d-1))\bigr).\end{align*}}%
    The last vertex is on $\ell$, and \[-\frac{1}{d-2}>-\frac{d}{d^2-d-1} > \cdots > -\frac{d^{m-2}}{d^{m-1}-d^{m-2}-1} > - \frac{d^{m-1}}{d^m-d^{m-1}-1}\] implies that all slopes in $N_d^-(\cF_0)$ is bigger than the slope of $\ell$, thus all other vertices are above $\ell$.
    
    \item $0<k<d-1$: We have \[\cF_k = (bd+1) \binom{d}{k} \sigma_m^k \tau_m^{d-1-k}.\] Since $\binom{d}{k}$ has $d$-valuation $1$, from \eqref{sigmaNewton} and (ii) the initial point of $N_d^-(\cF_k)$ is equal to $(m(d-1-k),d^m)$. Also, it consists of the line segments with slopes \[-\frac{1}{d-2} > - \frac{d}{d^2-d-1} > \cdots > -\frac{d^{m-2}}{d^{m-1}-d^{m-2}-1} > - \frac{D_{m-2}}{d^{m-2}}.\] Since the slope of $\ell$ is between the last two slopes, it suffices to show that the vertex between those two slopes is on or above $\ell$. From the initial point given above, this vertex is equal to \[(d^{m-1}(d-1-k),d^m-D_{m-1}(d-1-k)).\] It is on or above $\ell$ if and only if \[d^m-D_{m-1}(d-1-k) \ge -\frac{d^{m-1}}{d^m-d^{m-1}-1} (d^{m-1}(d-1-k) - (d^m-d^{m-1})) + d^m-d^{m-1},\] which is equivalent to \[d^m-d^{m-1}-1 \ge D_m k.\] We can check that this inequality is true for $k=d-2$, so for all $0<k<d-1$.
    
    \item $k=d-1$: We have \[\cF_{d-1} = bd(d-1)\sigma_m^{d-1},\] so using \eqref{sigmaNewton} we get \[N_d^- (\cF_{d-1}) = L\bigl((1,d^m),(d^{m-1}(d-1)+1,(d^{m-1}+1)(d-1)+1)\bigr).\] The first vertex is on $\ell$, and since \[-\frac{d^{m-1}}{d^m-d^{m-1}-1} > -\frac{D_{m-2}}{d^{m-2}}\] the other vertex is above $\ell$.
\end{enumerate}

Now it suffices to show that two vertices given in \eqref{eq:bdg} are actually vertices of $N_d^-(bd\cG_m)$, i.e., \[v_1(bd\cG_m) = d^m, \ \ v_{d^m-d^{m-1}} (bd\cG_m) = d^m-d^{m-1}.\] (A) and (B) says $b^2 \mid \cF_k$ for every $0\le k<d-1$, so (C) gives the first vertex. On the other hand, \eqref{sigmaNewton} says $V(\sigma_m) = d^{m-1}+1$ and (ii) says $V(\tau_m) = d^{m-1}$, so from (B) and (C) \[V(\cF_k) = 1 + k V(\sigma_m) + (d-1-k) V(\tau_m) \ge 1+(d-1)d^{m-1} = d^m-d^{m-1}+1\] for every $0<k\le d-1$. Therefore (A) gives the second vertex. This ends the proof.
\end{proof}
\end{theorem}

We are now ready to prove the main result.

\begin{proof}[Proof of Theorem 1]
From \eqref{gNewton}, there exists $F_m$, an irreducible factor of $G_m$ over $\bQ$ (or $\bQ_d$), whose degree is at least $d^m-d^{m-1}-1$. On the other hand, \eqref{gsigmatau} implies that \[\deg G_m \le (d-1)\max(\deg \sigma_m,\deg \tau_m).\] Therefore, \[\liminf_{m\to \infty} \frac{\deg F_m}{\deg G_m} \ge \liminf_{m\to \infty} \frac{(d^m-d^{m-1}-1)(d^d-1)}{d^{m+d-1}(d-1)}=\frac{d^d-1}{d^d}=1-\frac{1}{d^d}.\]

Now we prove the irreducibility of $G_m$ for $m\le d$. Since $m=2$ case is proved in \cite{han}, we can assume that $m\ge 3$. Then \[\deg \tau_m = \frac{d^{m+d-1}-d^{m-1}}{d^d-1}+1 = d^{m-1}+1, \ \ \deg \tau_{m+1} = d^m+1\] from Proposition \ref{sigmataudegree}, and consequently \[\deg \cG_m = d^m- d^{m-1}-1\] from Proposition \ref{taug}. Therefore, $G_m$ should be a constant multiple of $F_m$, i.e., $F_m$ is irreducible over $\bQ$.
\end{proof}

\end{document}